\documentclass{article} 
\usepackage{amsmath, amsthm , amsfonts, amssymb}
\usepackage{epsfig} 
\usepackage{epstopdf} 
\usepackage[colorlinks=true]{hyperref}
\hypersetup{urlcolor=blue, citecolor=red}
\allowdisplaybreaks

\newtheorem{theorem}{Theorem}[section]

\newtheorem{lemma}[theorem]{Lemma}
\newtheorem{proposition}{Proposition}

\theoremstyle{definition}

\usepackage{amssymb, amsthm, amsmath}

\title{KNOT AS A COMPLETE INVARIANT OF A MORSE-SMALE 3-DIFFEOMORPHISM WITH FOUR FIXED POINTS}

\begin{document}
\maketitle

\centerline{\scshape Olga Pochinka and Elena Talanova and Danila Shubin}
\medskip
{\footnotesize
 \centerline{National Research University}
   \centerline{Higher School of Economics}
   \centerline{Russia}
} 

\medskip

\bigskip


\begin{abstract}
Lens spaces are the only 3-manifolds that admit gradient-like flows with four fixed points. This is an immediate corollary of Morse inequality and of the Morse function with four critical points existence. A similar question for gradient-like diffeomorphisms is open. Solution can be approached by describing a complete topological conjugacy invariant of the class of considered diffeomorphisms and constructing of representative diffeomorphism for every conjugacy class by the abstract invariant. Ch. Bonnati and V. Z. Grines proved that the topological conjugacy class of Morse-Smale flows with unique saddle is defined by the equivalence class of the Hopf knot in $\mathbb S^2\times\mathbb S^1$ which is projection of one-dimensional saddle separatrice and used the mentioned approach to prove that the ambient manifold of a diffeomorphism of this class is the three-dimensional sphere. In the present paper similar result is obtained for the gradient-like diffeomorphisms with exactly two saddle points and the unique heteroclinic curve.
\end{abstract}


	\section{Introduction and formulation of results.}
Recall that a diffeomorphism $f\colon M^n\to M^n$, defined on an orientable connected closed smooth $n$-dimensional ($n\geq 1$) manifold $M^n$ is called {\it Morse-Smale diffeomorphism} if:
\begin{enumerate}
	\item its non-wandering set $\Omega_f$ consist of a finite number of hyperbolic orbits;
	\item intersection of the invariant manifolds $W^s_p$ and $W^u_q$ is transversal for any non-wandering points $p, q$.
\end{enumerate}
If $\sigma_1, \sigma_2$ are different saddle periodic points of the Morse-Smale diffeomorphism and $W^{s}_{\sigma_{1}}\cap W^{u}_{\sigma_{2}}\neq\emptyset$, then the intersection $W^{s}_{\sigma_{1}}\cap W^{u}_{\sigma_{2}}$ is called {\it heteroclinic intersection} and its connected components of dimension one are called {\it heteroclinic curves}.

In paper \cite{BGP} a complete topological classification of Morse-Smale diffeomorphisms on closed 3-manifolds was given, however its significant part is description of the invariant since it is designed for the wider class of diffeomorphisms. In some cases there are other more natural invariants that can be found without considering them as special case of the general one. Thus, in this paper, we establish that complete invariant for the wide class of 3-diffeomorphisms whose non-wandering set consists of four points is equivalence class of a Hopf knot in $\mathbb S^2\times\mathbb S^1$.

Recall, that {\it a knot} in manifold $\mathbb S^2\times\mathbb  S^1$ is a smooth embedding $\gamma\colon \mathbb S^1\to \mathbb S^2\times\mathbb S^1$ or the image $L=\gamma(\mathbb S^1)$ of the embedding. Knots $\gamma,\gamma'$ are called {\it smoothly homotopic} if there exists a smooth map $\Gamma\colon \mathbb S^1\times[0,1]\to \mathbb S^2\times\mathbb S^1$ such that $\Gamma(s,0)=\gamma(s)$ and $\Gamma(s,1)=\gamma'(s)$ for any $s\in\mathbb S^1$. If at the same time $\Gamma|_{\mathbb S^1\times\{t\}}$ is an embedding for any $t\in[0,1]$ then the knots are called {\it isotopic}. The knots $L,L'$ are called {\it equivalent} if there exists a homeomorphism $h\colon \mathbb S^2\times\mathbb S^1\to \mathbb S^2\times\mathbb S^1$ such that $h(L)=L'$. Let $[L]$ denote the equivalence class of $L$.

The knot $L\subset\mathbb S^2\times\mathbb S^1$ is {\it a Hopf knot} if the homomorphism $i_{L*}$ induced by inclusion $i_{L}\colon L\to\mathbb S^2\times\mathbb S^1$ is a group isomorphism for  $\pi_1(L)\cong\pi_1(\mathbb S^2\times\mathbb S^1)\cong\mathbb Z$. 

\begin{figure}[h!]
	\centerline
	{\includegraphics[width=10 cm]{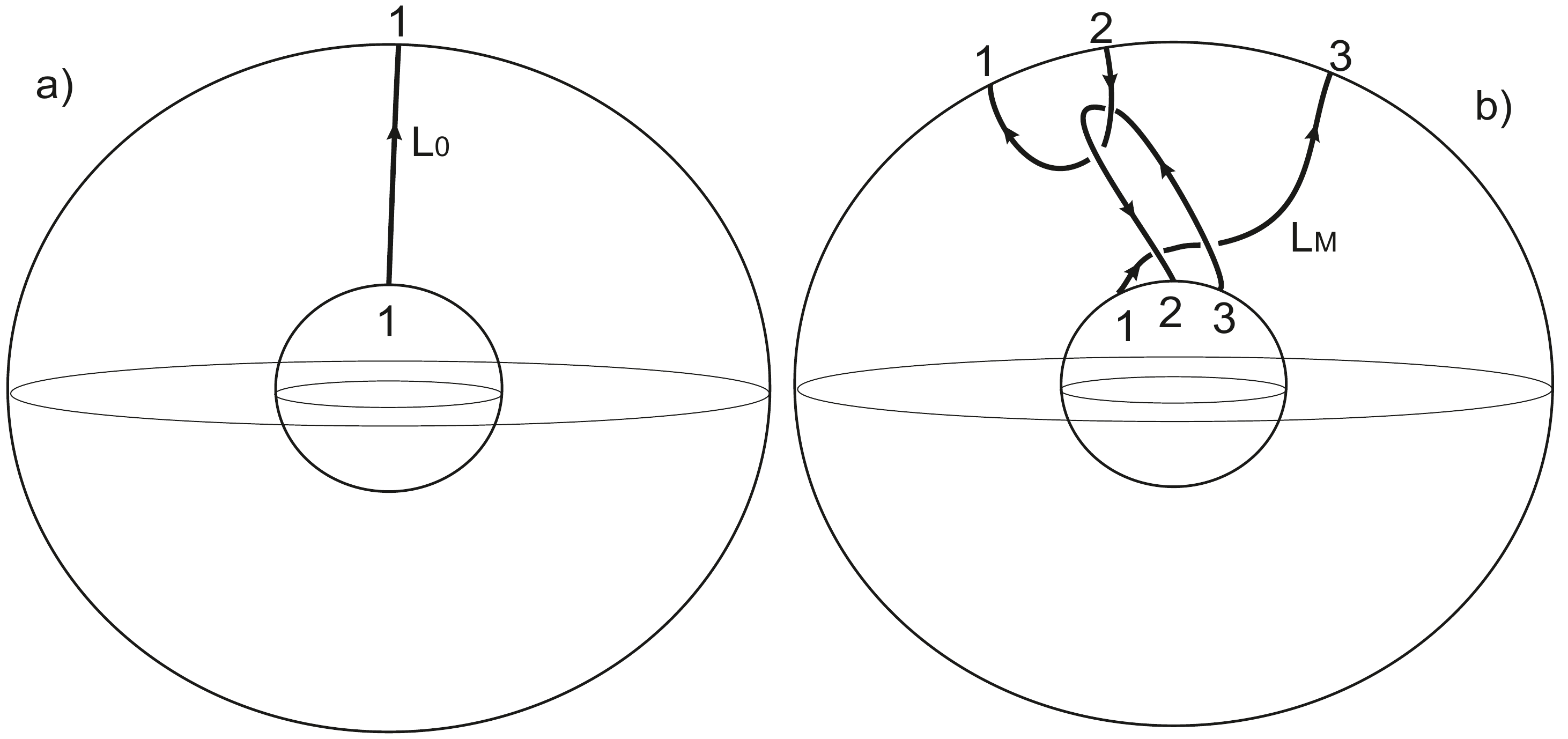}}
	\caption{Non-isotopic and non-equivalent Hopf knots $L_0$ and $L_M$: \\ a) standart Hopf knot $L_0$; b) Masur knot $L_{M}$.}\label{not}
\end{figure}

\begin{figure}[h!]
	\begin{minipage}[b]{0.45\textwidth}
		\includegraphics[width=\textwidth]{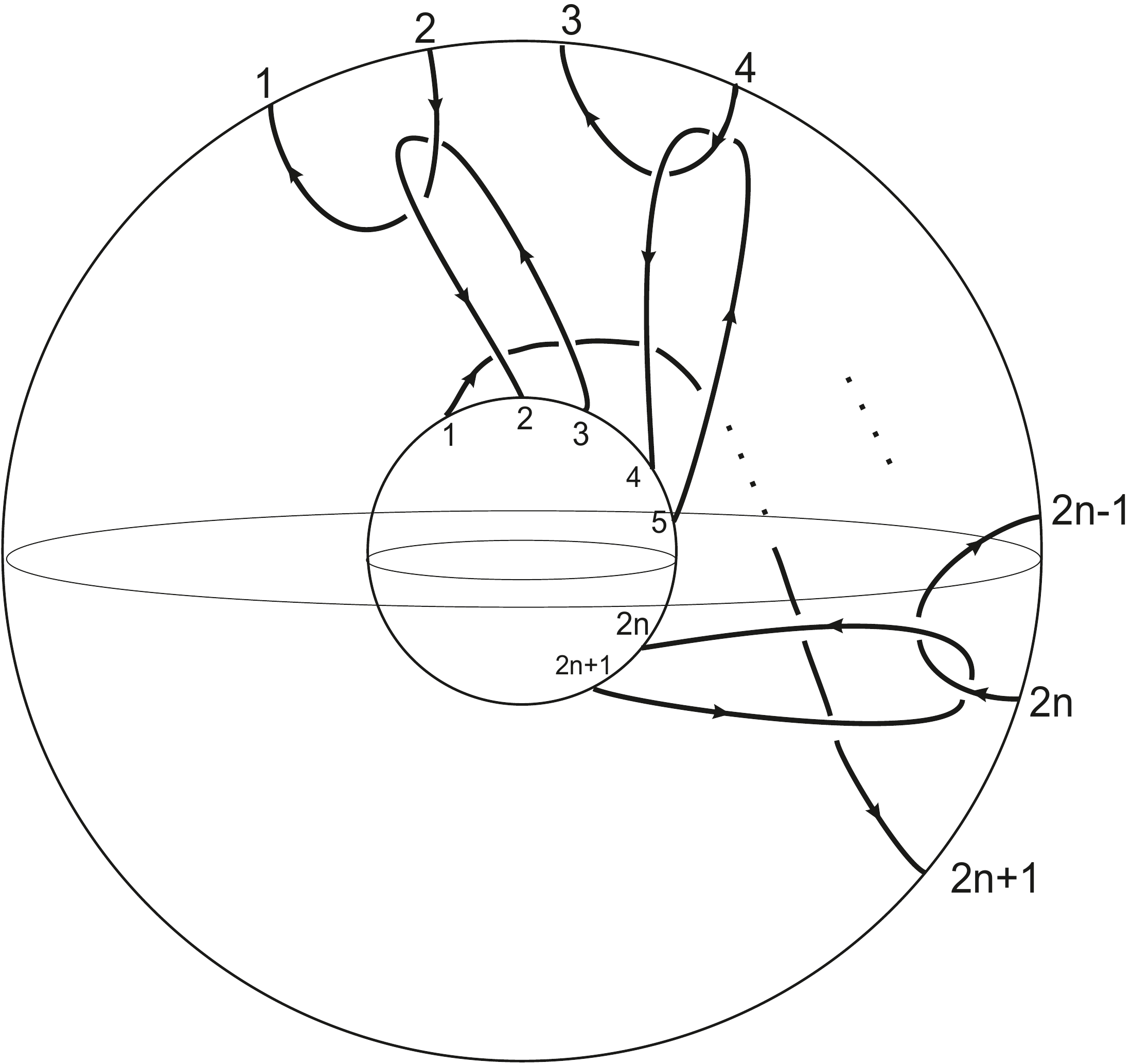}
		\center{Generalized Masur knot $L_{M,n}$}
	\end{minipage}
	\begin{minipage}[b]{0.45\textwidth}
		\includegraphics[width=\textwidth]{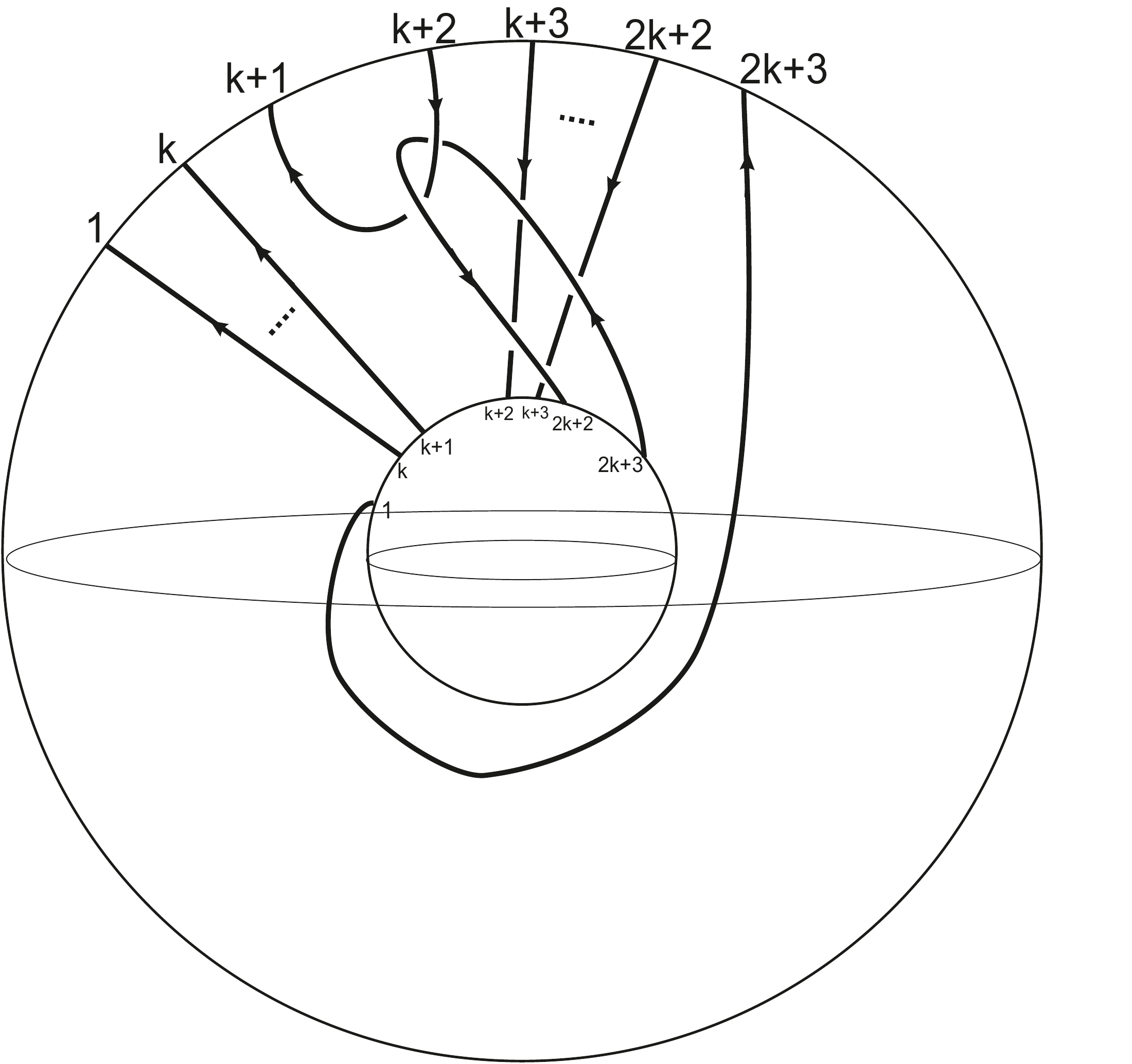}
		\center{Generalized Masur knot $L^k_{M}$}
	\end{minipage}
	\caption{A countable family of pairwise non-equivalent Hopf knots.}
	\label{Lmn}
\end{figure}

Note, that any Hopf knot is smoothly homotopic to the standard Hopf knot $L_0=\{x\}\times\mathbb S^1$ (see~e.g.~\cite{K-L}), but generally they are neither isotopic nor equivalent. It is known (see~e.g.~\cite{AMP}), that Hopf knots are equivalent if and only if they are isotopic. 
B.~Masur constructed a Hopf knot $L_M$ which is neither equivalent nor isotopic 
to $L_0$ (see~Pic.~\ref{not}).
In paper \cite{AMP} a countable  family of pairwise non-equivalent Hopf knots was constructed (see~Pic.~\ref{Lmn}).

Three-dimensional Morse-Smale diffeomorphisms with exactly four non-wandering points can be divided into two classes: in the first class each diffeomorphism has one saddle point, in the second class each diffeomorphism has two saddle points. It was established in \cite{BoGr} that topological conjugacy class of a diffeomorphism of the first class is defined by equivalence class of the Hopf knot, which is projection of one-dimensional saddle separatrix. According to \cite{Pi} and \cite{BoGr} any Hopf knot can be realized as a diffeomorphism on 3-sphere from the first class.

In the present paper diffeomorphisms of the second class are considered. Namely, the class $G$ of orientation-preserving Morse-Smale diffeomorphisms on closed manifold $M^3$ with the following properties:
\begin{itemize}
	\item non-wandering set of the diffeomorphism $f\in G$ consists of exactly four fixed points $\omega_f, \sigma_f^1, \sigma_f^2, \alpha_f$ and dimensions of their unstable manifolds are $0, 1, 2, 3$ respectively;
	\item the set $H_f=W^{s}_{\sigma_f^{1}}\cap{W^{u}_{\sigma_f^{2}}}$ is not empty and is path-connected (hence, it consists of the unique non-compact curve)\footnote{In the paper \cite{GMZ} it was proved that for any diffeomorphism $f$ of the second class the set $H_f$ contains at least one non-compact heteroclinic curve, but generally it contain more curves, including the case, when it contains  infinitely many of compact heteroclinic curves.}.
\end{itemize}

\begin{figure}[h!]
	\center{\includegraphics
		[width=0.65\linewidth]{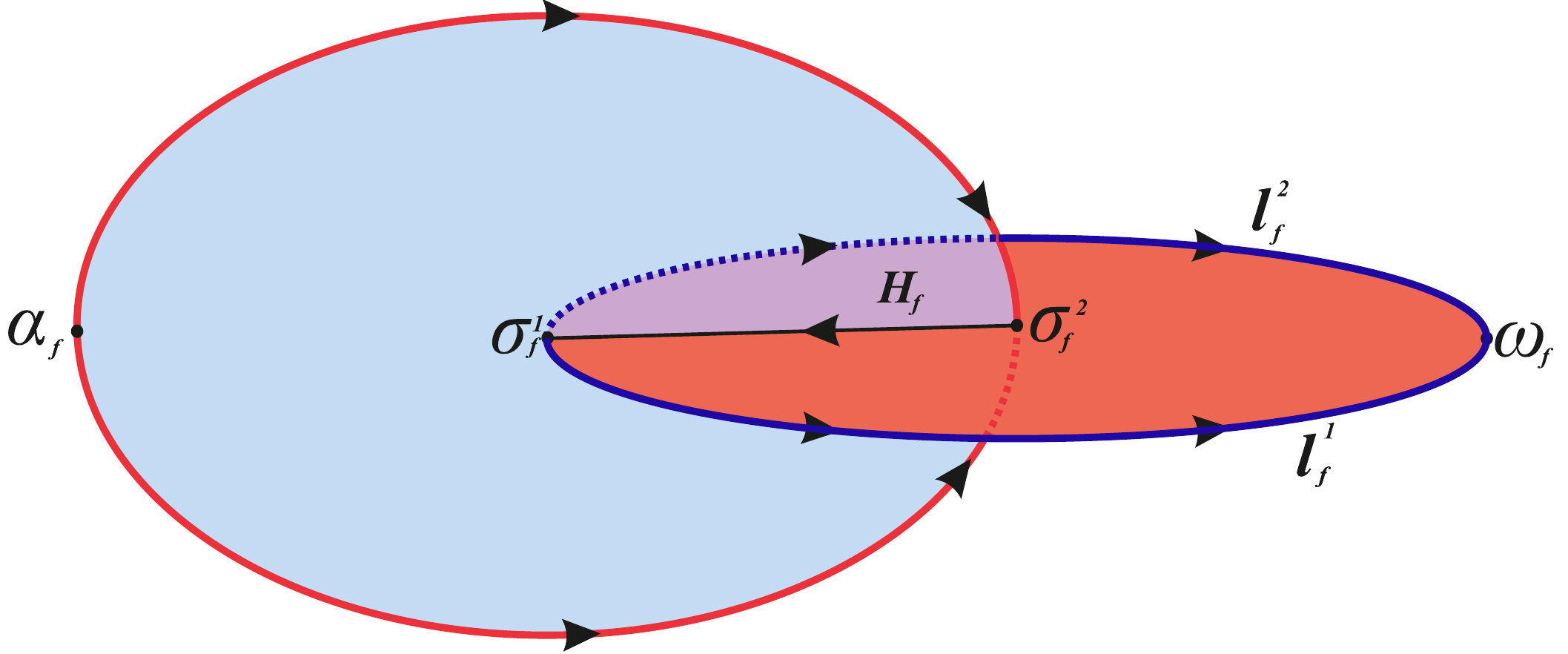}}
	\caption{Dynamics of the diffeomorphism $f\in G$}
	\label{pic1}
\end{figure}

Let $f\in G$. Denote the unstable separatrices of the point $\sigma_f^{1}$ as $\ell_f^1$, $\ell_f^2$. Thus (see~e.g.~\cite{GrMePo2016}), the closure $cl(\ell_f^i)$ ($i=1,2$) of one-dimensional unstable separatrix of the point $\sigma_f^{1}$ is homeomorphic to the simple closed curve and it consists of the separatrix and two points: the saddle $\sigma_f^{1}$ and the sink $\omega_f$ (see~Pic.~\ref{pic1}).

Let ${\bf x}=(x_1,x_2,x_3)\in\mathbb R^3$, $|| {\bf x} ||=\sqrt{x_1^2 + x_2^2 + x_3^2}$ and $a\colon\mathbb R^3\to\mathbb R^3$ be a diffeomorphism given by the formula $a({\bf x})=\frac{\bf x}{2}.$ Define the map $p\colon\mathbb R^3\setminus O\to\mathbb S^{2}\times\mathbb S^1$ by the formula $$p({\bf x})=\left(\frac{x_1}{||{\bf x}||},\frac{x_2}{||{\bf x}||}, 	\log_2(||{\bf x}||)\pmod 1\right).$$
Let $V_{\omega_f}=W^s_{\omega_f}\setminus\omega_f$.
By virtue of the hyperbolicity of the sink $\omega_f$, there exists a diffeomorphism $\psi_f\colon V_{\omega_f}\to\mathbb R^3\setminus O$, which conjugates the diffeomorphisms $f$ and $a$. Let $p_{\omega_f}=p\psi_f:V_{\omega_f}\to\mathbb S^{2}\times\mathbb S^1$ and $L_f^{i}=p_{\omega_f}(\ell_f^i),\,i=1,2$. 

\begin{lemma}\label{odin}
	For any diffeomorphism $f\in G$ the sets $L_f^1,\,L_f^2$ are equivalent Hopf knots in $\mathbb S^{2}\times\mathbb S^1$. 
\end{lemma}
Let $\mathcal L_f=[L_f^1]=[L_f^2]$ the equivalence class of this knots.

\begin{theorem}\label{iffknot}
	The diffeomorphisms $f,f'\in G$ are topologically conjugated if and only if $\mathcal L_f=\mathcal L_{f'}$.
\end{theorem}

So, the equivalence class of the Hopf knot in $\mathbb S^{2}\times\mathbb S^1$ is a complete topological invariant for  diffeomorphisms from $G$. Moreover, the following theorem holds.

\begin{theorem}\label{Pidy}
	For any equivalence class $\mathcal L$ of Hopf knot in $\mathbb S^{2}\times\mathbb S^1$ there exists a diffeomorphism  $f_{\mathcal L}:\mathbb S^3\to\mathbb S^3\in G$ such that $\mathcal L_{f_{\mathcal L}}=\mathcal L$.
\end{theorem}

The immediate corollary of Theorem~\ref{iffknot} and Theorem~\ref{Pidy} is the fact, that the ambient manifold of the diffeomorphism of the class $G$ is the 3-sphere $\mathbb S^3$ (an independent proof of this fact see in \cite{PoSh}).

\section{Compatible foliated neighborhoods.}\label{so}

For $t\in(0,1]$ let $\mathcal N_{1}^t=\{(x_1,x_2,x_3)\in\mathbb{R}^3~:~ x_1^2(x_2^2+x_3^2)< t\}$, $\mathcal N_{2}^t=\{(x_1,x_2,x_3)\in\mathbb{R}^3~:~ (x_1^2+x_2^2)x_3^2< t\}$ and for $i\in\{1,2\}$ let $\mathcal N^1_{i}=\mathcal N_{i}$.  

In the neighborhood $\mathcal N_{1}$ define a pair of transversal foliations $\mathcal{F}^u_1,~\mathcal{F}^s_{1}$ in the following way:

$$\mathcal{F}^u_1=\bigcup\limits_{(c_{2},c_3)\in Ox_2x_3}\{(x_1,x_2,x_3)\in \mathcal N_{1}~:~(x_{2},x_3)=(c_{2},c_3)\},$$ $$\mathcal{F}^s_{1}=\bigcup\limits_{c_1\in Ox_1}\{(x_1,x_2,x_3)\in \mathcal N_{1}~:~x_1=c_1\}.$$

In the neighborhood $\mathcal N_{2}$ define a pair of transversal foliations $\mathcal{F}^u_2,~\mathcal{F}^s_{2}$ in the following way: 
$$\mathcal{F}^u_2=\bigcup\limits_{c_3\in Ox_3}\{(x_1,x_2,x_3)\in \mathcal N_{2}~:~x_3=c_3\},$$ $$\mathcal{F}^s_{2}=\bigcup\limits_{(c_{1},c_2)\in Ox_1x_2}\{(x_1,x_2,x_3)\in \mathcal N_{2}~:~(x_{1},x_2)=(c_{1},c_2)\}.$$

Define the diffeomorphisms $a_i\colon\mathbb R^3\to\mathbb R^3$ by the formula: $$a_1({\bf x})=\left(2x_1,\frac{x_2}{2},\frac{x_3}{2}\right),\,a_2=a_1^{-1}.$$
Note, that for $i\in\{1,2\}$ the set $\mathcal N_{i}^t$ is invariant under the diffeomorphism $a_i$ action, which maps the leaves of $\mathcal{F}^u_i$ ($\mathcal{F}^s_{i}$) to leaves of $\mathcal{F}^u_i$ ($\mathcal{F}^s_{i}$).

By virtue of \cite{BGP}, the saddle point $\sigma_f^i$ of the diffeomorphism $f\in G$ has an {\it linearizing neighborhood} $N_f^{i}$, which is equipped with a homeomorphism ${\mu}_{i}\colon N_{i}\to {\mathcal N}_{i}$ conjugating the diffeomorphism $f\vert_{{N}_f^{i}}$ with the diffeomorphism $a_i|_{{\mathcal N}_{i}}$ and being a diffeomorphism on $N_i\setminus(W^s_{\sigma_f^i}\cup W^u_{\sigma_f^i})$. The foliations $\mathcal{F}^u_{i},~\mathcal{F}^s_{i}$ induce $f$-invariant foliations ${F}^u_{i},~{F}^s_{i}$ on the linearizing neighborhood $N_{i}$ by homeomorphism ${\mu}_{i}^{-1}$. Let ${F}^u_{i,x}$ (${F}^s_{i,x}$) denote the unique leave of the foliation ${F}^u_{i}$ (${F}^s_{i}$) containing the point $x\in N_{i}$. 

Moreover, the heteroclinic curve $H_f$ possesses an $f$-invariant neigborhood $N_{H_f}$ with two-dimensional foliation $F$, whose every leaf $F_x$ containing a point $x\in N_{H_f}$ transversally intersect $H_f$ at the unique point. Also the foliation $F$ is transversal to the foliations  $F^s_1,\,F^u_2$ simultaneously. 

\begin{figure}[h!]
	\centerline
	{\includegraphics[width=12 cm]{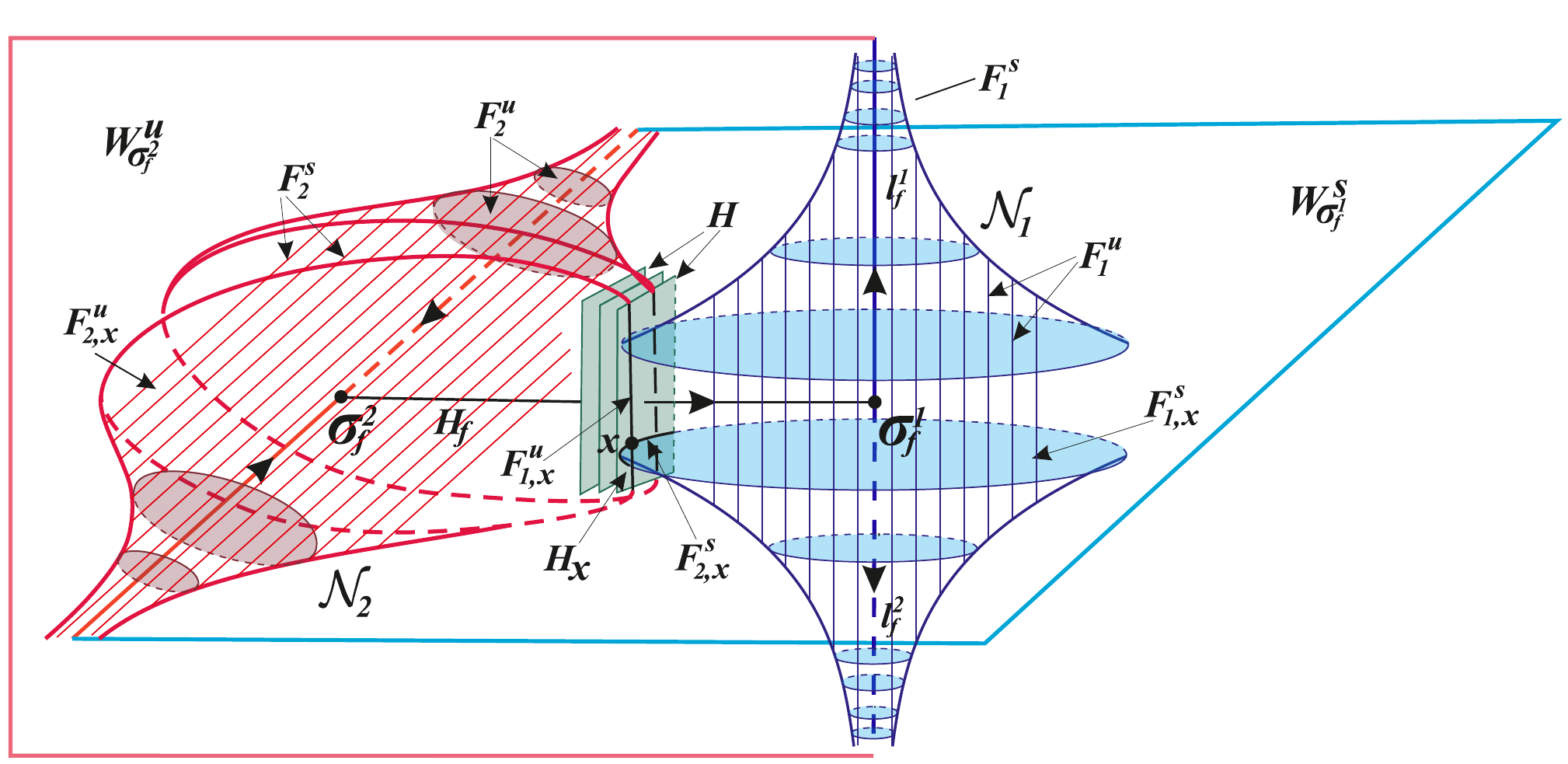}}
	\caption{Compatible foliated neighborhoods.}\label{sogl}
\end{figure}

The neighborhoods $N_{H_f},\,N_1,\,N_2$ with foliations $F,\,F^s_1,\,F^u_2$ 
are called {\it compatible foliated neighborhoods} if for any point $x\in(N_f^1\cap N_f^2\cap N_{H_f})$ and the leave $F_x$ of the foliation $F$ 
the following conditions hold (see Pic. \ref{sogl}):
$${F}^s_{1,x}\cap F_x={F}^s_{2,x}\cap({N}_f^{1}\cap N_{H_f}),\quad {F}^u_{2,x}\cap F_x={F}^u_{1,x}\cap
({N}_f^{2}\cap N_{H_f}).$$
\begin{proposition}[\cite{BGP}, Theorem 1] For any diffeomorphism $f\in G$ there exists compatible foliated neighborhoods. 
\end{proposition}

\section{Equivalence of the knots $L^1_f,\,L^2_f$.}

In this section we prove Lemma~\ref{odin}:  for any diffeomorphism $f\in G$ the sets $L_f^1,\,L_f^2$ are equivalent Hopf knots in $\mathbb S^{2}\times\mathbb S^1$.

\begin{proof}
As it was mentioned in the introduction, Hopf knots are equivalent if and only if they are isotopic. So, to prove the lemma it is sufficient to construct an isotopy between the knots $L_f^1,\,L_f^2$. Let (see~Pic.~\ref{Cf}): $${C}_f =p_{\omega_f}(W^u_{\sigma_f^{2}}).$$ 

	\begin{figure}[h!]
		\centerline
		{\includegraphics[width=12 cm]{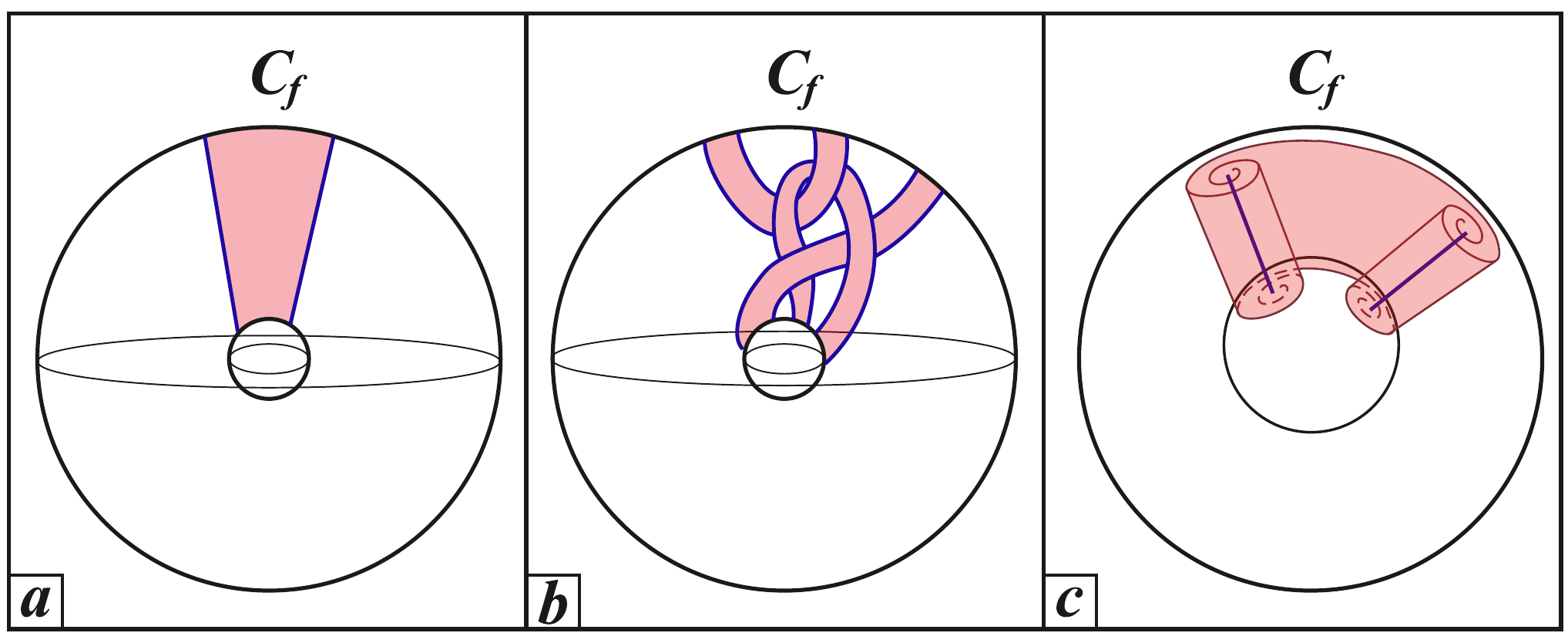}}
		\caption{Possible variants of the projection $C_f$ of the two-dimensional unstable saddle manifold}\label{Cf}
	\end{figure}

Since the diffeomorphism $f|_{W^u_{\sigma_f^{2}}}$ is topologically conjugated with the linear extension, the orbit space  $(W^u_{\sigma_f^{2}}\setminus\sigma_f^{2})/f$ is homeomorphic to two-dimensional torus. 
Since $W^u_{\sigma_f^{2}}\cap V_{\omega_f}=W^u_{\sigma_f^{2}}\setminus (H_f\cup \sigma_f^{2})$ and orbit space $H_f/f$ is homeomorphic to the circle, the set $C_f$ is homeomorphic to two-dimensional annulus. 
Moreover, the homomorphism $i_{C_f*}$, induced by inclusion $i_{C_f}\colon C_f\to\mathbb S^2\times\mathbb S^1$ is group isomorphism of $\pi_1(C_f)\cong\pi_1(\mathbb S^2\times\mathbb S^1)\cong\mathbb Z$. Let $$U^1_f=p_{\omega_f}(N^1_f).$$ 
Since $N^1_f\cap V_{\omega_f}=N^1_f\setminus W^s_{\sigma_f^{1}}$, the set $U^1_f$ is disjoint union of two solid tori $U^1_f=U^{1,1}_f\sqcup U^{1,2}_f$ which are tubular neighborhood of the knots $L^1_f,\,L^2_f$. 
Let (see~Pic.~\ref{Ti}) $$T^i_f=\partial U^{1,i}_f,\,S^i_f=T^i_f\cap C_f.$$ 

	\begin{figure}[h!]
		\centerline
		{\includegraphics[width=5 cm]{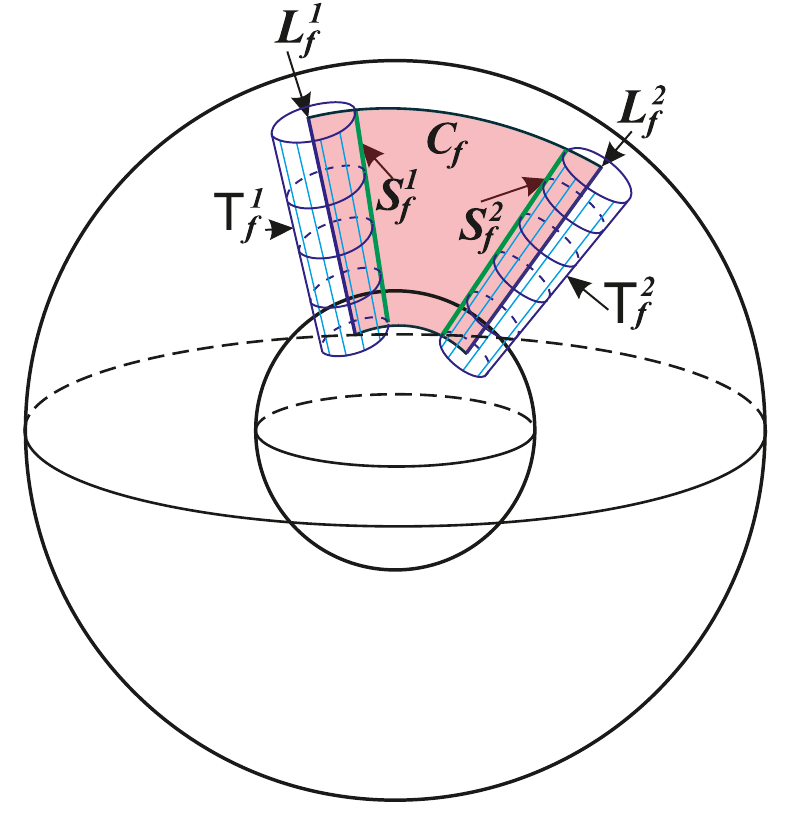}}
		\caption{The projection of compatible neighborhood $N_1$}\label{Ti}
	\end{figure}

	Since the set $W^{s}_{\sigma_f^{1}}\cap{W^{u}_{\sigma_f^{2}}}$ consists of the unique non-compact $f$-invariant curve, the set $\partial N^{1}_{f}\cap{W^{u}_{\sigma_f^{2}}}$ consists of two non-compact $f$-invariant curves as well. The projections of this curves in $\mathbb S^2\times\mathbb S^1$ are curves $S^1_f\sqcup S^2_f$. It implies that $S^1_f,\, S^2_f$ are isotopic Hopf knots in $\mathbb S^2\times\mathbb S^1$. 
	
	So, the knots $S^i_f,\,L^i_f$ are the generators of the solid torus $U^{1,i}_f$. Thus, they bound a two-dimensional annulus in the solid torus and consequently are isotopic.
\end{proof}

\section{Hopf knot equivalence class as a complete invariant of the topological conjugacy in class $G$.}

In this section we  prove Theorem~\ref{iffknot}: the diffeomorphisms $f,f'\in G$ are topologically conjugated if and only if $\mathcal L_f=\mathcal L_{f'}$.
\begin{proof} $ $

	$\Rightarrow$ Let the diffeomorphisms $f,f'\in G$ be topologically conjugated by a homeomorphism $h\colon M^3\to M^3$. Since $h$ maps the invariant manifolds of the diffeomorphism $f$  fixed points into invariant manifolds of the diffeomorphism $f'$ fixed points preserving their stability, we get $h(W^s_{\omega_f}) = W^s_{\omega_{f'}}, h(\ell^i_f) = \ell^i_{f'}$. Since $hf=f'h$, the homeomorphism $h$ defines the homeomorphism $\hat h\colon\mathbb S^2\times\mathbb S^1\to\mathbb S^2\times\mathbb S^1$ by the formula:
	$$\hat h=p_{\omega_f}hp_{\omega_f}^{-1}.$$
	It implies that $\hat h(L^i_f)=L^i_{f'}$ and consequently the Hopf knots $L^i_f,\,L^i_{f'}$ are equivalent.
	
	$\Leftarrow$ Let $\mathcal L_f=\mathcal L_{f'}$. Thus there exists a homeomorphism $\hat h_0\colon \mathbb S^2\times\mathbb S^1\to\mathbb S^2\times\mathbb S^1$ such that $\hat h_0(L^1_f)=L^1_{f'}$. We will modify the homeomorphism $\hat h_0$ step by step to construct a homeomorphism $h\colon M^3\to M^3$ conjugating the diffeomorphisms $f$ and $f'$. Doing that we will use the notation of Lemma~\ref{odin} and the section~\ref{so} which will be equipped with prime mark for the diffeomorphism $f'$. Let (see~Pic.~\ref{V}):
	$$U^2_f=p_{\omega_f}(N^2_f),\,U_f=U^1_f\cup U^2_f.$$  

	{\bf Step 1. Construction of the homeomorphism $\hat h_1:\mathbb S^2\times\mathbb S^1\to\mathbb S^2\times\mathbb S^1$ such that $\hat h_1(U_f)=U_{f'}$.}
By virtue of Lemma~\ref{odin}, $U_f$ is tubular neighborhood of the knot $L^1_f$ ($L^2_f$ also). Let $\tilde U_f\supset U_f$  be a  tubular neighborhood of the knot $L^1_f$ also. Let $U=\hat h_0^{-1}(U_{f'})$, $\tilde U=\hat h_0^{-1}(\tilde U_{f'})$ and chose a  tubular neighborhood $V$ of the knot $L^1_f$ such that $V\subset int(U_f\cap U)$ (see~Pic.~\ref{V}).
	\begin{figure}[h!]
		\centerline
		{\includegraphics[width=5 cm]{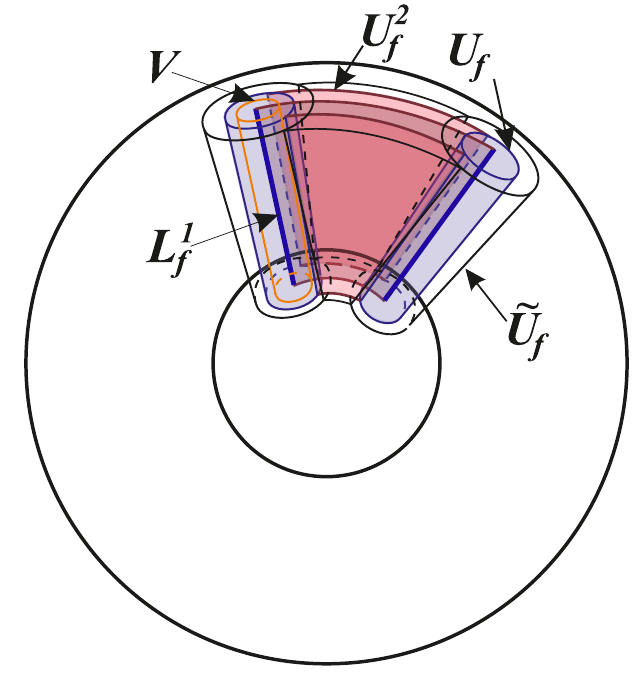}}
		\caption{Tubular neighborhoods of the knot $L^1_f$}\label{V}
	\end{figure}
	
	Since the sets $\tilde U_f\setminus int\,U_f$ and $U_f\setminus int\,V$ are homeomorphic to $\mathbb T^2\times[0,1]$, then there is a homeomorphism $\hat\psi_f:\mathbb S^2\times\mathbb S^1\to\mathbb S^2\times\mathbb S^1$, which coincides with identity outside $\tilde U_f$ and such that $\hat\psi_f(U_f)=V$. In the same way construct the homeomorphism $\hat\psi\colon \mathbb S^2\times\mathbb S^1\to\mathbb S^2\times\mathbb S^1$, which coincides with identity outside $\tilde U$ and such that $\hat\psi(U)=V$. 	Thus, the desired homeomorphism $\hat h_1$ is $$\hat h_1=\hat h_0\hat\psi^{-1}\hat\psi_f.$$
	
	{\bf Step 2. Construction of the homeomorphism $\hat h_2:\mathbb S^2\times\mathbb S^1\to\mathbb S^2\times\mathbb S^1$, which coincide with $\hat h_1$ outside $\hat U_f$ and such that $\hat h_2(C_f)=C_{f'}$.}
	
	Let $H=\mu_1(H_f),\,H'=\mu'_1(H_{f'})$. Thus, $H,\,H'$ are $a_1$-invariant curves on the plane $Ox_2x_3$. Thus, there exists a homeomorphism $\xi_s\colon Ox_2x_3\to Ox_2x_3$, which commutes with $a_1|_{Ox_2x_3}$ and such that $\xi_s(H)=H'$ (see~e.g.~\cite{GrMePo2016}). Define the homeomorphism $\xi\colon\mathbb R^3\to\mathbb R^3$ by the formula (see~Pic.~\ref{xi}):
	$$\xi(x_1,x_2,x_3)=(x_1,\xi_s(x_2,x_3)).$$ 
	\begin{figure}[h!]
		\centerline
		{\includegraphics[width=10 cm]{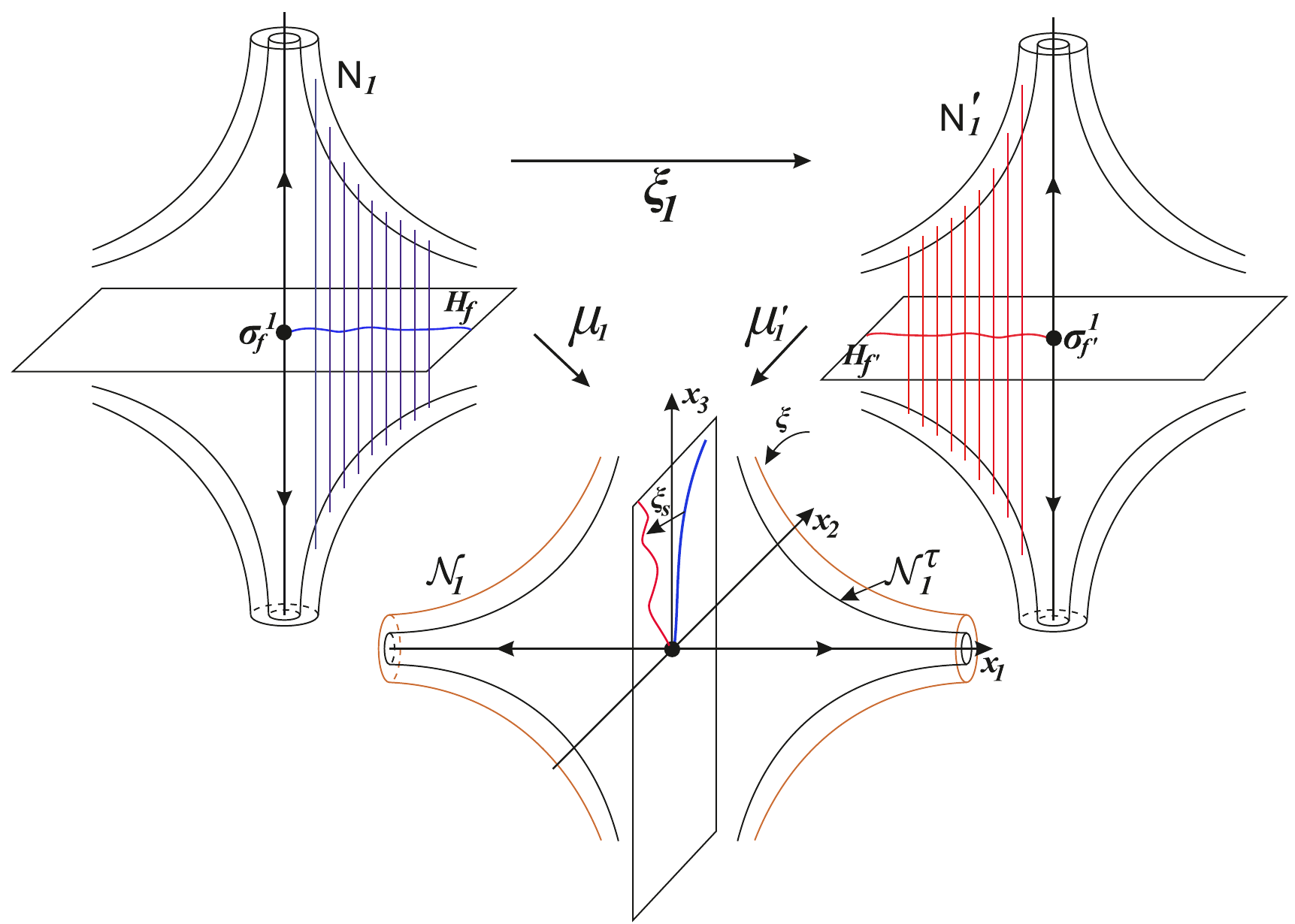}}
		\caption{Construction of the homeomorphism $\xi$}\label{xi}
	\end{figure}
	
	Chose $\tau\in(0,1)$ to get $\xi(\mathcal N^\tau_1)\subset int\,\mathcal N_1$. Let $N^\tau_1=\mu_1(\mathcal N^\tau_1)$ and $\xi_1 = \mu'^{-1}_1 \xi \mu_1|_{N^\tau_1}$. The definition of compatible system of neighborhoods implies that $\xi_1(N^\tau_1\cap W^u_{\sigma^2_f})\subset W^u_{\sigma^2_{f'}}$. Let $N^\tau_2=\mu_2(\mathcal N^\tau_2)$, $U^{\tau}_1=p_{\omega_f}(N^\tau_1),\,U^{\tau}_2=p_{\omega_f}(N^\tau_2),\,U^\tau=U^\tau_1\cup U^\tau_2,\,U'^{\tau}_2=p_{\omega_{f'}}(N^\tau_2)$ and $\hat\xi_1=p_{\omega_{f'}}\xi_1p_{\omega_f}^{-1}|_{U^\tau_1}$. Let $U'^\tau_1=\hat\xi_1(U^\tau_1),\,U'^\tau=U'^\tau_1\cup U'^{\tau}_2$, $K=C_f\setminus int\,U^\tau,\,K'=C_{f'}\setminus int\,U^\tau_1,\, U_K=U^\tau\setminus U^\tau_1,\,U_{K'}=U'^\tau\setminus U'^\tau_1$ (see~Pic.~\ref{xi1})
	
	\begin{figure}[h!]
		\centerline
		{\includegraphics[width=8 cm]{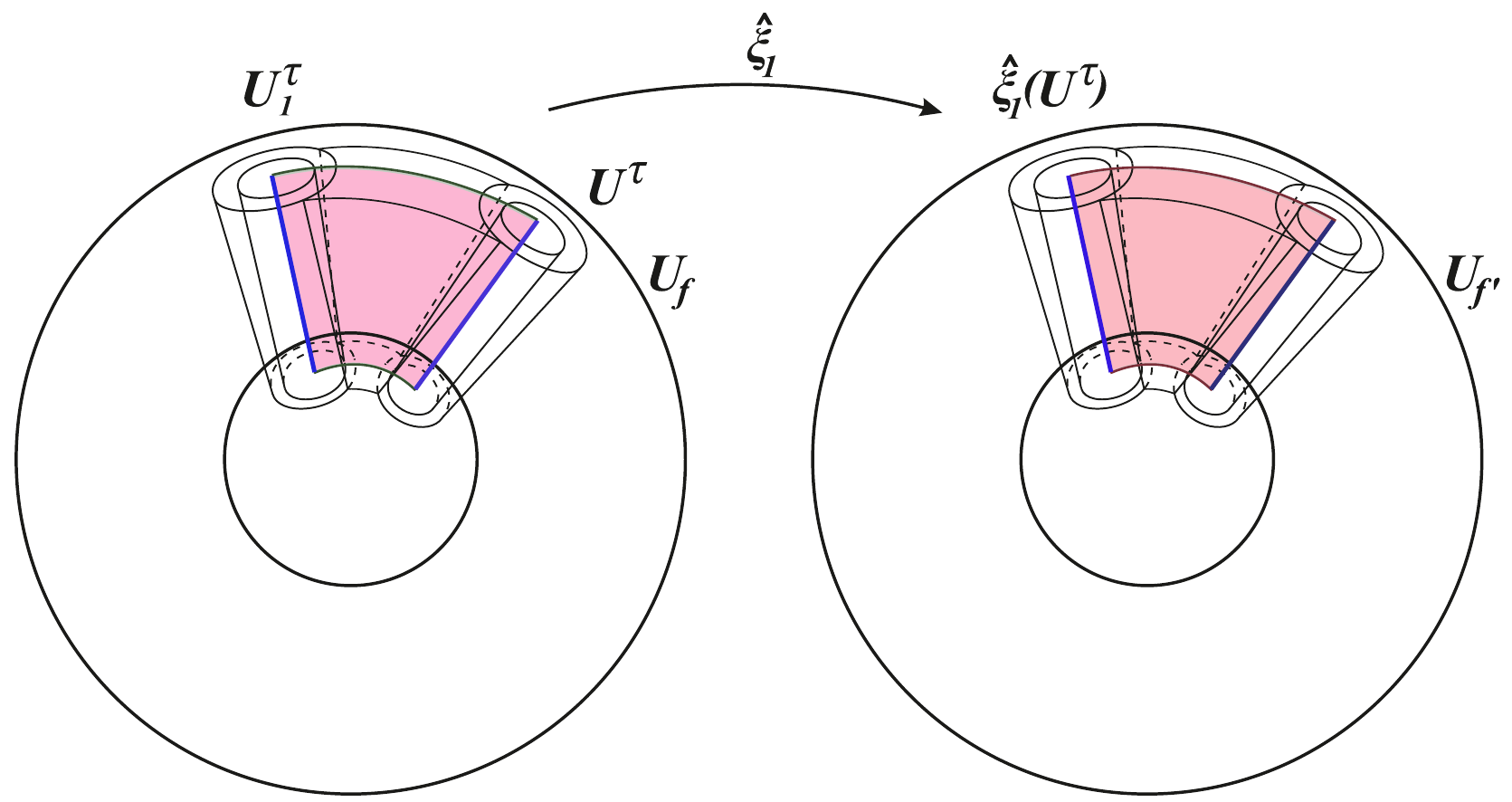}}
		\caption{Construction of the homeomorphism $\hat\xi_1$}\label{xi1}
	\end{figure}
	
	Since the sets $K,\,K'$ are homeomorphic to two-dimensional annuli and the sets $U_K,\,U_{K'}$ are their tubular neighborhoods, the homeomorphism $\hat\xi_1$ can be extended to the homeomorphism $\hat\xi\colon U^\tau\to U'^\tau$ such that $\hat\xi(C_f)=C_{f'}$. 
	Thus, the homeomorphism $\hat\psi_1=\hat h_1^{-1}\hat\xi$ possess the property $\hat\psi_1(U^\tau)\subset int\, U_f$ and the homeomorphism $\hat\xi_1|_{\partial U^\tau}$ is homotopic to identity. 
	Since the sets $U_f\setminus int\,U^\tau$ and $U_f\setminus int\,\hat\psi_1(U^\tau)$ are homeomorphic to $\mathbb T^2\times[0,1]$, 
	the homeomorphism $\hat\psi_1$ extends to the homeomorphism $\hat\psi_1\colon\mathbb S^2\times\mathbb S^1\to\mathbb S^2\times\mathbb S^1$, which is identity outside $U_f$. 
	Thus, the desired homeomorphism $\hat h_2$ is $$\hat h_2=\hat h_1\hat\psi_1.$$ 
	
	{\bf Step 3. Construction of the desired homeomorphism $h$.}
	By the construction of $\hat h_2$ there exists its lift $h_2\colon V_{\omega_f}\to V_{\omega_{f'}}$, conjugating the diffeomorphism $h_2\colon V_{\omega_f}\to V_{\omega_{f'}}$ and extending to the homeomorphism $\xi_1$ on $W^s_{\sigma_1}$. So, the conjugating homeomorphism is defined everywhere except the closures of the one-dimensional manifolds of the saddle points. By virtue of the \cite[Theorem 1]{BGP} such homeomorphism can be extended to the desired homeomorphism $h$.
\end{proof}

\section{Realization of the diffeomorphisms of the class $G$.}
The Fox-Artin arc first encountered in dynamics in the paper of D. Pixton~\cite{Pi}. In that paper the Morse-Smale diffeomorphism on 3-sphere with unique saddle point, whose invariant manifolds form the Fox-Artin arc, was constructed. In the paper~\cite{BoGr} arbitrary Hopf knot in $\mathbb S^{2}\times\mathbb S^1$ was realized as Morse-Smale diffeomorphism with unique saddle point on the 3-sphere (see~e.g.~\cite{GrMePo2016} and~\cite{MePo2018}). In the present section we give similar realization for diffeomorphisms of the class $G$.

\begin{figure}[h!]
	\centerline{\includegraphics[width=5 true cm, height=7.5 true cm]{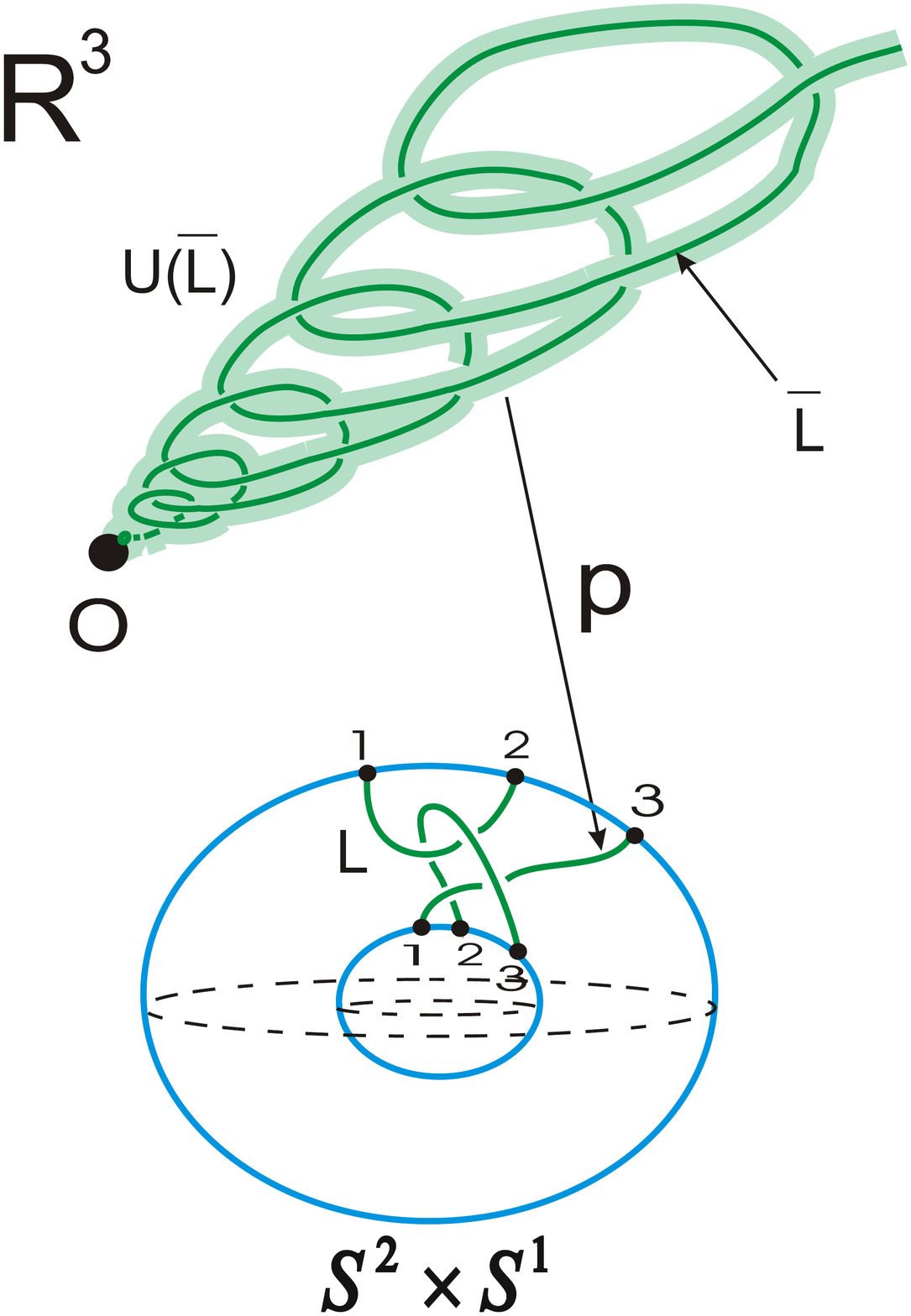}}
	\caption{Hopf knot lift}\label{af}
\end{figure}

Let $L\subset \mathbb S^{2}\times\mathbb S^1$ be Hopf knot and  $U_L$ be its tubular neighborhood. Thus, the set $\bar L=p^{-1}(L)$ is $a$-invariant curve in $\mathbb R^3$ and $U_{\bar L}=p^{-1}(U_L)$ its $a$-invariant neighborhood, diffeomorphic to $\mathbb{D}^{2}\times\mathbb R^1$ (see~Pic.~\ref{af}).
Define the flow $g^t\colon C\to C$ on the cylinder $C=\{{\bf x}\in\mathbb R^3~:~x_2^2+x_{3}^2\leqslant 4\}$  by the formula $$g^t({\bf x})=(x_1+t,x_2,x_3).$$
Thus, there exists a diffeomorphism ${\zeta}\colon {U_L}\to C$ which conjugates $a\vert_{{U_L}}$ and $g=g^1|_C$. {Define the flow $\phi^t$ on $C$ by the formula
	$$\begin{cases}
		\dot{x}_1=\begin{cases}1-\frac{1}{9}(||{\bf x}||-4)^2, \quad ||{\bf x}|| \leqslant 4 \cr
			1, \quad ||{\bf x}|| > 4
		\end{cases}\cr
		\dot{x}_2=\begin{cases}
			\frac{x_2}{2}\big(\sin\big(\frac{\pi}{2}\big(||{\bf x}||-3\big)\big)-1\big), \quad 2<||{\bf x}||\leqslant 4\cr
			-x_2,\quad \quad ||{\bf x}||\leqslant 2\cr
			0, \quad ||{\bf x}|| > 4
		\end{cases}\cr
		\dot{x}_3=\begin{cases}
			-\frac{x_3}{2}\big(\sin\big(\frac{\pi}{2}\big(||{\bf x}||-3\big)\big)-1\big), \quad 2<||{\bf x}||\leqslant 4\cr
			x_3,\quad \quad ||{\bf x}||\leqslant 2\cr
			0, \quad ||{\bf x}|| > 4.
		\end{cases}
	\end{cases}$$} 
\begin{figure}[h!]
	\centerline{\includegraphics[width=8 true cm, height=7.5 true cm]{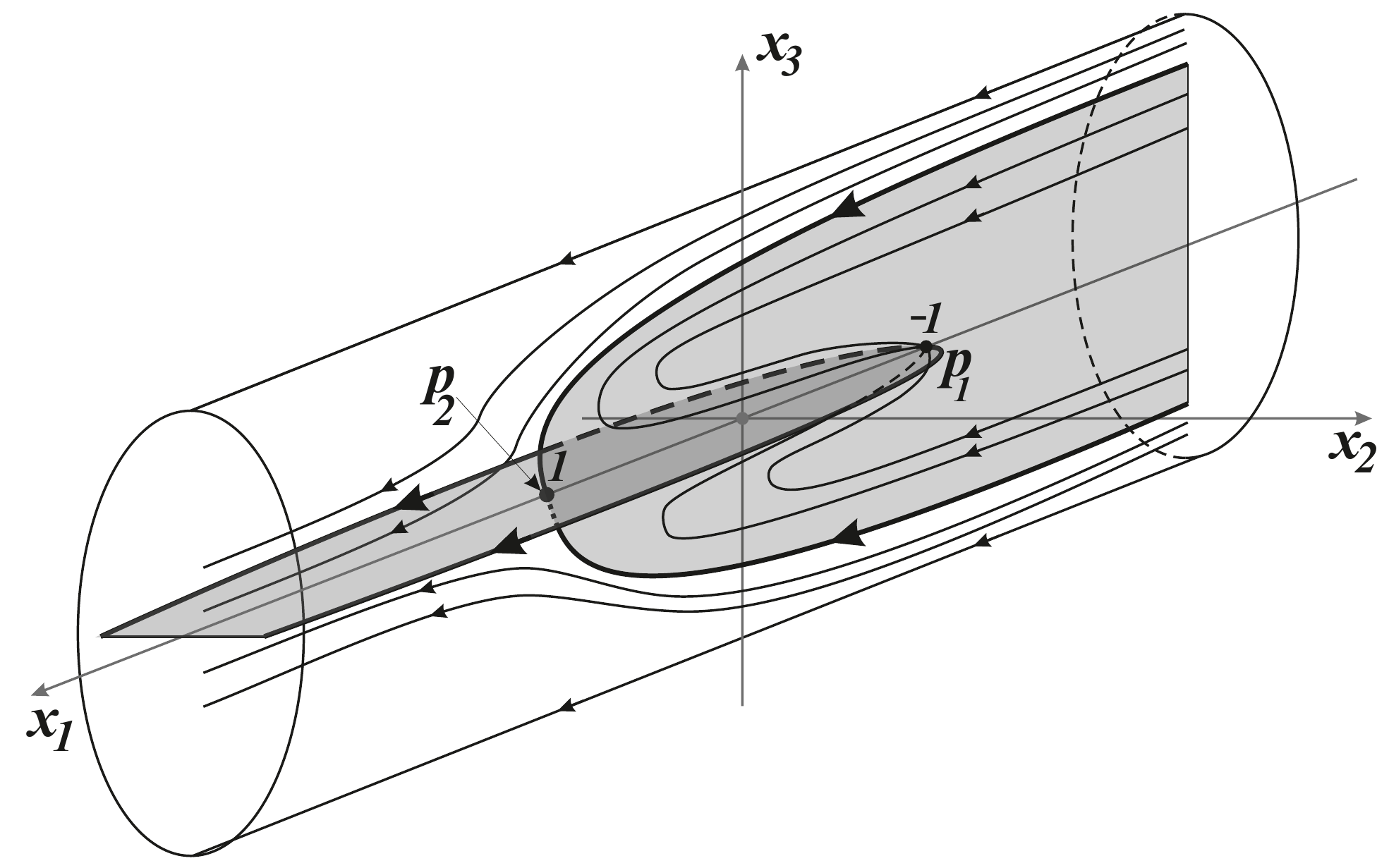}}
	\caption{Trajectories of the flow $\phi^t$}\label{cherry}
\end{figure}

Construction of the diffeomorphism $\phi=\phi^1$ implies that it has two hyperbolic fixed saddles: point $P_1(-1,0,0)$ with Morse index $1$ and the point $P_2(1,0,0)$ with Morse index $2$ (see~Pic.~\ref{cherry}). The non-compact heteroclinic curve $W^{s}_{P_1}\cap W^{u}_{P_2}$ coincides with the open interval $\left\{{\bf x}\in\mathbb R^3:\,|x_1|<1,\,x_2=x_3=0\right\}$. Note, that $\phi$ coincides with the diffeomorphism $g=g^1$ outside the ball $\{{\bf x}\in C:||{\bf x}||\leqslant 4\}$.

\begin{figure}[h!]
	\centerline{\includegraphics[width=6 true cm, height=6 true cm]{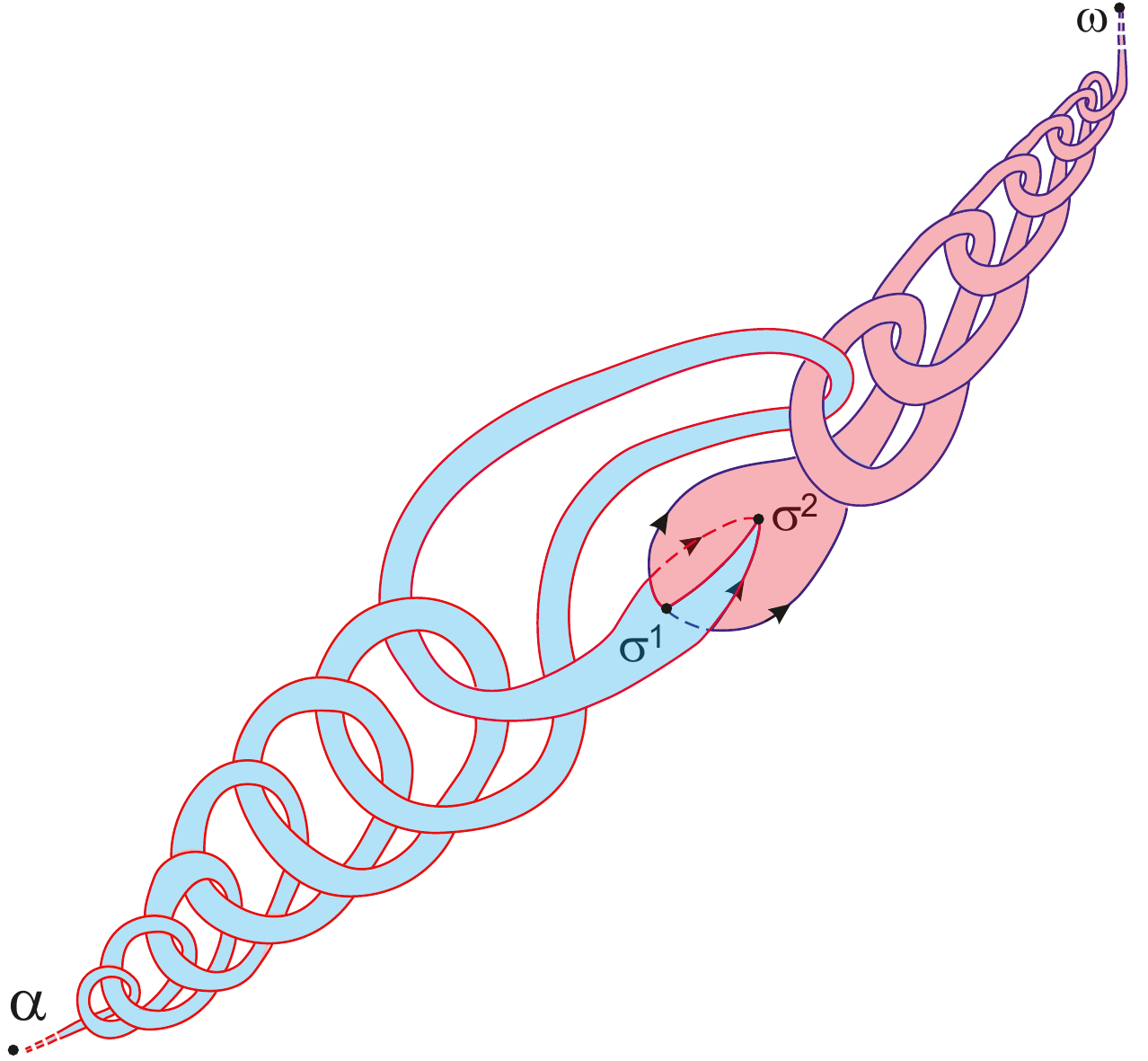}}
	\caption{Phase portrait of the diffeomorphism $\bar f_L$}\label{af2}
\end{figure}

Define the diffeomorphism $\bar f_L\colon\mathbb R^3\to\mathbb R^3$ such that $\bar{f}_{L}$ coincides with $a$ outside ${U_{L}}$ and coincides with ${\zeta}^{-1}\phi{\zeta}$ on ${U_{L}}$. Thus, $\bar f_{L}$ has two hyperbolic fixed points in ${U_{L}}$: the saddle ${\zeta}^{-1}(P_1)$ and the saddle ${\zeta}^{-1}(P_2)$ (see~Pic.~\ref{af2}).

Denote the North pole of the sphere $\mathbb S^3$ with $N(0,0, 0, 1)$ and the standard stereographic projection with $\vartheta:\mathbb R^3\to(\mathbb{S}^3\setminus\{N\})$.
{The construction of the diffeomorphism $\bar{f}_{L}$ implies that it coincides with the diffeomorphism $a$ in some neighborhood of point $O$ and outside some neighborhood of this point. Thus, it induces the Morse-Smale diffeomorphism $f_{{L}}$ on $\mathbb{S}^3$: $$f_{{L}}(s)=\begin{cases}\vartheta(\bar f_{L}(\vartheta^{-1}(s))),~{s}\neq N;\cr N,~{s}=N\end{cases}.$$}
It follows directly from the construction that non-wandering set of the diffeomorphism $f_{{L}}$ consists of four hyperbolic fixed points: one sink $\omega$, two saddles $\sigma^1=\vartheta({\zeta}^{-1}(P_1))$, $\sigma^2=\vartheta({\zeta}^{-1}(P_2))$ and one source $\alpha$. The diffeomorphism belongs to the class $G$ and $\mathcal L_{f_L}=[L]$.

\section*{Acknowledgments} 

The study was founded by Russian Scientific Foundation, agreement \textnumero 22-11-00027 except the realization of the considered class diffeomorphisms, which was founded by Laboratory of dynamical systems and applications NRU HSE, of the Ministry of science and higher education of the RF grant ag. \textnumero 075-15-2019-1931. Particular thanks to our colleague Andrey Morozov for the beautiful pictures.



\begin{thebibliography}{1}
	
	\bibitem{BGP}
	C.~Bonatti, V.~Grines, and O.~Pochinka, ``Topological classification of
	morse-smale diffeomorphisms on 3-manifolds,'' {\em Duke Mathematical
		Journal}, vol.~168, no.~13, pp.~2507--2558, 2019.
	
	\bibitem{K-L}
	P.~Kirk and C.~Livingston, ``Knot invariants in 3-manifolds and essential
	tori,'' {\em Pacific Journal of Mathematics}, vol.~197, no.~1, pp.~73--96,
	2001.
	
	\bibitem{AMP}
	P.~Akhmet’ev, T.~Medvedev, and O.~Pochinka, ``On the number of the classes of
	topological conjugacy of pixton diffeomorphisms,'' {\em Qualitative Theory of
		Dynamical Systems}, vol.~20, no.~3, pp.~1--15, 2021.
	
	\bibitem{BoGr}
	C.~Bonatti and V.~Grines, ``Knots as topological invariants for gradient-like
	diffeomorphisms of the sphere {$S^3$},'' {\em Journal of Dynamical and
		Control Systems}, vol.~6, no.~4, pp.~579--602, 2000.
	
	\bibitem{Pi}
	D.~Pixton, ``Wild unstable manifolds,'' {\em Topology}, vol.~16, pp.~167--172,
	12 1977.
	
	\bibitem{GMZ}
	V.~Z. Grines, E.~V. Zhuzhoma, and V.~S. Medvedev, ``On morse--smale
	diffeomorphisms with four periodic points on closed orientable manifolds,''
	{\em Mathematical Notes}, vol.~74, no.~3, pp.~352--366, 2003.
	
	\bibitem{GrMePo2016}
	V.~Grines, T.~Medvedev, and O.~Pochinka, {\em Dynamical Systems on 2- and
		3-Manifolds}, vol.~46.
	\newblock 01 2016.
	
	\bibitem{PoSh}
	V.~I. Shmukler and O.~V. Pochinka, ``Bifurcations that change the type of
	heteroclinic curves of the morse-smale 3-diffeomorphism,'' {\em Taurida
		Journal of Computer Science Theory and Mathematics}, no.~1 (50),
	pp.~101--114, 2021.
	
	\bibitem{MePo2018}
	T.~V. Medvedev and O.~V. Pochinka, ``The wild fox-artin arc in invariant sets
	of dynamical systems,'' {\em Dynamical Systems}, vol.~33, no.~4,
	pp.~660--666, 2018.
	
\end{thebibliography}
\end{document}